\documentclass{amsart}
\usepackage{latexsym}
\usepackage{amsfonts}
\usepackage{amssymb}
\usepackage{enumerate}
\usepackage{amssymb,latexsym,amsxtra,amscd}
\usepackage{amsmath,amsthm,amsfonts, amssymb,amscd}
\usepackage[all]{xy}

\newtheorem{theorem}{Theorem}[section]
\newtheorem{lemma}[theorem]{Lemma}%[section]
\newtheorem{definition}[theorem]{Definition}%[section]
\newtheorem{example}[theorem]{Example}%[section]
\newtheorem{remark}[theorem]{Remark}%[section]
\newcommand{\Ker}{{\rm Ker}}
\newcommand{\Id}{{\rm Id}}
\newcommand{\Rad}{{\rm Rad}}
\newcommand{\MV}{{\rm MV}}

\long\def\alert#1{\smallskip{\hskip\parindent\vrule%
\vbox{\advance\hsize-2\parindent\hrule\smallskip\parindent.4\parindent%
\narrower\noindent#1\smallskip\hrule}\vrule\hfill}\smallskip}

\begin{document}
\title[State-Morphism Algebras - General Approach]{State-Morphism Algebras -
General Approach}
\author[M. Botur and A. Dvure\v{c}enskij]{Michal Botur$^1$ and Anatolij Dvure\v{c}enskij$^2$}
\date{}
\maketitle

\begin{center}  \footnote{


MB thanks for the support by SAIA, Slovakia,  by MSM 6198959214 of
the RDC of the Czech Government, and by GA\v{C}R P201/11/P346, Czech
Republic,   AD thanks  for the support by Center of Excellence SAS
-~Quantum Technologies~-,  ERDF OP R\&D Projects CE QUTE ITMS
26240120009 and meta-QUTE ITMS 26240120022, the grant VEGA No.
2/0032/09 SAV. }
\small{Department of Algebra and Geometry\\
Faculty of  Science, Palack\'y University\\
17. listopadu 12, CZ-771 46 Olomouc, Czech Republic\\
$^2$ Mathematical Institute,  Slovak Academy of Sciences\\
\v Stef\'anikova 49, SK-814 73 Bratislava, Slovakia\\

E-mail: {\tt michal.botur@upol.cz},\ {\tt dvurecen@mat.savba.sk} }
\end{center}

\begin{abstract} We present a complete description of subdirectly
irreducible state BL-algebras as well as of subdirectly irreducible
state-morphism BL-algebras. In addition, we present a general theory
of state-morphism algebras, that is, algebras of general type with
state-morphism which is an idempotent endomorphism. We define a
diagonal state-morphism algebra and we show that every subdirectly
irreducible state-morphism algebra can be embedded into a diagonal
one. We describe generators of varieties of state-morphism algebras,
in particular ones of state-morphism BL-algebras, state-morphism
MTL-algebras, state-morphism non-associative BL-algebras, and
state-morphism pseudo MV-algebras.

\end{abstract}

\section{Introduction}%1

A state, as an analogue of a probability measure, is a basic notion
of the theory of quantum structures, see e.g. \cite{DvPu}. However,
for MV-algebras, the state as averaging the truth value in the \L
ukasiewicz logic was introduced firstly by Mundici in \cite{Mun}, 40
years after introducing MV-algebras, \cite{Cha}. We recall that a
state on an MV-algebra $\mathbf M$ is a mapping $s:M\to [0,1]$ such
that (i) $s(a\oplus b) = s(a)+s(b),$ if $a\odot b =0,$ and (ii)
$s(1)=1.$ The property (i) says that $s$ is additive on mutually
excluding events $a$ and $b.$ It is important note that every
non-degenerate MV-algebra admits at least one state. The set of
states is a convex set, which in the weak topology of states is a
compact Hausdorff set, and every extremal state is in fact an
MV-algebra homomorphism from $\mathbf M$ into the MV-algebra of the
real interval $[0,1],$ and vice-versa, \cite{Mun}. In addition,
extremal states generate the set of all states because by the
Krein-Mil'man Theorem, \cite[Thm 5.17]{Goo}, every state is a weak
limit of a net of convex combinations of these special
homomorphisms.

In the last decade, the states entered into theory of MV-algebras in
a very ambitious manner. In \cite{Mun1, KuMu}, authors have showed
a relation between states and  de Finetti's approach to probability
in terms of bets. In addition, Panti and independently Kroupa in
\cite{Pan, Kro} have  showed that every state on $\mathbf M$ is an
integral through a unique regular Borel probability measure
concentrated on the set of extremal states on $\mathbf M$.

Nevertheless as we have seen states are not a proper notion of
universal algebra, and therefore, they do not provide  an
algebraizable logic for probabilistic reasoning of the many-valued
approach.

Recently, Flaminio and Montagna in \cite{FlMo} presented an
algebraizable logic containing probabilistic  reasoning, and its
equivalent algebraic semantic is the variety of state MV-algebras.
We recall that a {\it state MV-algebra} is an MV-algebra whose
language is extended adding an operator, $\tau$ (called also an {\it
internal state}),  whose properties are inspired by the ones of
states. The analogues of extremal states are {\it state-morphism
operators}, introduced in \cite{DiDv}. By definition, it is an
idempotent endomorphism on an MV-algebra.

State MV-algebras generalize, for example, H\'ajek's approach,
\cite{Haj}, to fuzzy logic with modality Pr (interpreted as {\it
probably}) which has the following semantic interpretation: The
probability of an event $a$ is presented as the truth value of
Pr$(a)$. On the other hand, if $s$ is a state, then $s(a)$ is
interpreted as averaging of appearing the many valued event $a.$

We note that if $(\mathbf M,\tau)$ is a state MV-algebra, assuming
that that the range $\tau(\mathbf M)$ is simple, we see that it is a
subalgebra of the real interval $[0,1]$ and therefore,  $\tau$ can
be regarded as a standard state on $\mathbf M.$ On the other hand,
every MV-algebra $\mathbf M$ can be embedded into the tensor product
$[0,1]\otimes \mathbf  M,$ therefore, given a state $s$ on $\mathbf
M,$ we define an operator $\tau_s$ on $[0,1]\otimes \mathbf M$ via
$\tau_s(t\otimes a) := t\cdot s(a),$ \cite[Thm 5.3]{FlMo}.  Then due
to \cite[Thm 3.2]{DiDv}, $\tau_s$ is a state-operator that is a
state-morphism operator iff $s$ is an extremal state. Thus, there is
a natural correspondence between the notion of  a state and an
extremal state on one side, and a state-operator and a
state-morphism operator on the other side.

Subdirectly irreducible state-morphism MV-algebras were described in
\cite{DiDv, DDL2} and this was extended also for state-morphism
BL-algebras in \cite{Dvu1}. A complete description of both
subdirectly irreducible state MV-algebras as well as subdirectly
irreducible state-morphism MV-algebras can be found in \cite{DKM}.
In \cite{DDL1}, it was shown that if $(\mathbf M,\tau)$ is a state
MV-algebra whose image $\tau(\mathbf M)$ belongs to the variety
generated by the $L_1,\ldots, L_n,$ where $L_i :=\{0,1/i,\ldots,
i/i\},$ then $\tau$ has to be a state-morphism operator. The same is
true if $\mathbf M$ is linearly ordered, \cite{DiDv}.  Recently, in
\cite{DKM}, we have shown that the unit square $[0,1]^2$ with the
diagonal operator generates the whole variety of state-morphism
MV-algebras; it answered in positive an open problem posed in
\cite{DiDv}. In addition, there was shown that in contrast to
MV-algebras, the lattice of subvarieties is uncountable. Moreover,
it was shown that every subdirectly irreducible state-morphism
MV-algebra can be embedded into some diagonal one.

In this  paper, we continue in the study of state BL-algebras and
state-morphism BL-algebras. Because the methods developed in
\cite{DKM} are so general that, it is possible to study more general
structures than MV-algebras or BL-algebras under a common umbrella.
Hence, we introduce state-morphism algebras $(\mathbf A,\tau),$
where the algebra $\mathbf A$ is an arbitrary algebra of type $F$
and $\tau$ is an idempotent endomorphism of $\mathbf A.$ Then
general results applied to special types of algebras give
interesting new results.

The main goals of the paper are:

(1) Complete characterizations of subdirectly irreducible state
BL-algebras and state-morphism BL-algebras.

(2) Showing that every subdirectly state-morphism algebra can be
embedded into some diagonal one $D(\mathbf B):=(\mathbf B \times
\mathbf B,\tau_B),$ where $\tau(a,b)=(a,a),$ $a,b \in B,$ which is
also subdirectly irreducible.

(3) We show that if $\mathcal K$ is a generator of some variety
$\mathcal V$ of algebras of type $F,$ then the system of diagonal
state-morphism algebras $\{D(\mathbf B): \mathbf B \in \mathcal K\}$
is a generator of the variety of state-morphism algebras whose
$F$-reduct belongs to $\mathcal V.$

(4)  We exhibit cases when the Congruence Extension Property holds
for a variety of state-morphism algebras.

(5) In particular, a generator of the variety of state-morphism
BL-algebras is the class of all BL-algebras of the real interval
$[0,1]$ equipped with a continuous t-norm. Similarly, a generator of
the variety of state-morphism MTL-algebras is the class of all
MTL-algebras of the real interval equipped with a left-continuous
t-norm, similarly for non-associative BL-algebras one is the set of
all non-associative BL-algebras of the real interval $[0,1]$
equipped with a non-associative t-norm, and a generator of the
variety of state-morphism pseudo MV-algebras is any pseudo
MV-algebra $\Gamma(G,u),$ where $(G,u)$ is a doubly transitive
unital $\ell$-group.

\section{Subdirectly Irreducible State BL-algebras}%2

In this section, we define state BL-algebras and state-morphism
BL-algebras and we present a complete description of their
subdirectly irreducible algebras. These results generalize those
from \cite{DiDv, DDL2, Dvu1, DKM}.

We recall that according to \cite{Haj}, a {\it BL-algebra} is an
algebra $\mathbf M=(M; \wedge, \vee, \odot, \rightarrow, 0,1)$ of
the type $\langle 2,2,2,2,0,0\rangle$ such that $(M;
\wedge,\vee,0,1)$ is a bounded lattice, $(M;\odot,1)$ is a
commutative monoid, and for all $a,b,c \in M,$

\begin{enumerate}
\item[(1)] $c\leq a\rightarrow b$ iff $a\odot c\leq b;$
\item[(2)] $a\wedge b=a\odot (a\rightarrow b);$
\item[(3)] $(a\rightarrow b)\vee (b\rightarrow a)=1.$
\end{enumerate}

For any $a \in M,$ we define a complement $a^-:=a\to 0.$  Then $a
\le a^{--}$ for any $a \in M$ and a BL-algebra is an MV-algebra iff
$a^{--} =a$ for any $a \in M.$

A non-empty set $F\subseteq M$ is called a {\it filter} of $
\mathbf{M}$ (or a {\it BL-filter} of $\mathbf{M}$) if for every
$x,y\in M$: (1) $x,y\in F$ implies $x\odot y\in F,$ and (2) $x\in
F,$ $x\leq y$ implies $y\in F.$ A  filter $F \ne M$ is called a {\it
maximal filter} if it is not strictly contained in any other filter
$F'\ne M.$  A BL-algebra is called {\it local} if it has a unique
maximal filter.

We denote by $\Rad_1(\mathbf{M})$ the intersection of all maximal
filters of $\mathbf{M}.$

Let $\mathbf{M}$ be a BL-algebra. A mapping $\tau:M \to M$ such
that, for all $x,y \in M,$ we have

\begin{enumerate}

\item[$(1)_{BL}$] $\tau(0)=0;$
\item[$(2)_{BL}$] $\tau(x\rightarrow y)=\tau(x)\rightarrow \tau(x\wedge y);$
\item[$(3)_{BL}$] $\tau(x\odot y)=\tau(x)\odot \tau(x\rightarrow (x\odot y));$
\item[$(4)_{BL}$] $\tau(\tau(x)\odot \tau(y))=\tau(x)\odot \tau(y);$
\item[$(5)_{BL}$] $\tau(\tau(x)\rightarrow \tau(y))=\tau(x)\rightarrow
\tau(y)$
\end{enumerate}
is said to be a {\it state-operator} on $\mathbf{M},$ and the pair
$(\mathbf{M},\tau)$ is said to be a {\it state BL-algebra}, or more
precisely, a {\it BL-algebra with  internal state}.

If $\tau:M\to M$ is a BL-endomorphism such that $\tau\circ \tau =
\tau,$ then $\tau$ is a state-operator on $\mathbf{M}$ and it is
said to be a {\it state-morphism operator} and the couple
$(\mathbf{M},\tau)$ is said to be a {\it state-morphism BL-algebra}.

A filter $F$ of a BL-algebra $\mathbf{M}$ is said to be a
$\tau$-filter if $\tau(F)\subseteq F.$ If $\tau$ is a state-operator
on $\mathbf{M},$ we denote by
$$
\Ker(\tau)=\{a\in M: \tau(a)=1\}.
$$
then $\Ker(\tau)$ is a $\tau$-filter. A state-operator  $\tau$ is
said to be {\it faithful} if $\Ker(\tau) =\{1\}.$

We recall that there is a one-to-one relation  between congruences
and $\tau$-filters on a state BL-algebra $(\mathbf{M},\tau)$ as
follows. If $F$ is a $\tau$-filter, then the relation $\sim_F$ given
by $x \sim_F y$ iff $x\to y,y\to x\in F$ is a congruence of the
BL-algebra $\mathbf{M}$ and $\sim_F$ is also a congruence of the
state BL-algebra $(\mathbf{M},\tau).$

Conversely, let $\sim$ be a congruence of state BL-algebra
$(\mathbf{M},\tau)$ and set $F_\sim:=\{x\in M : x\sim 1\}.$ Then
$F_\sim$ is a $\tau$-filter of $(\mathbf{M},\tau)$ and
$\sim_{F_\sim} = \sim$ and $F= F_{\sim_F}.$

By \cite[Lem 3.5(k)]{CDH}, $(\tau(\mathbf M),\tau)$ is a subalgebra
of $(\mathbf{M},\tau),$ $\tau$ on $\tau(M)$ is the identity, and
hence, $(\Ker(\tau); \to, 0,1)$ is a subhoop of $\mathbf{M}$. We say
that two subhoops, $A$ and $B,$ of a BL-algebra $\mathbf{M}$ have
the \emph{disjunction property} if for all $x\in A$ and $y\in B$, if
$x\vee y=1$, then either $x=1$ or $y=1$.

Nevertheless a subdirectly irreducible state BL-algebra
$(\mathbf{M},\tau)$ is not necessarily linearly ordered, according
to \cite[Thm 5.5]{CDH}, $\tau(\mathbf M)$ is always linearly
ordered.

We note that according to  \cite[Prop 3.13]{CDH}, if $\mathbf{M}$ is
an MV-algebra, then the notion of a state MV-algebra coincides with
the notion of a state BL-algebra.

The following three characterizations were originally proved in
\cite{DKM} for state MV-algebras. Here we show that the original
proofs from \cite{DKM} slightly improved work also for state
BL-algebras.

\begin{lemma}\label{le:2.1}
Suppose that $(\mathbf{M},\tau )$ is a  state BL-algebra. Then:
\begin{enumerate}
\item[{\rm (1)}] If $\tau $ is faithful, then $(\mathbf{M},\tau)$ is a
subdirectly irreducible state BL-algebra if and only if $\tau
(\mathbf{M})$ is a subdirectly irreducible BL-algebra.
\end{enumerate}

Now let $(\mathbf{M},\tau)$ be subdirectly irreducible. Then:

\begin{enumerate}
\item[{\rm (2)}] $\Ker(\tau)$ is (either trivial or) a subdirectly
irreducible hoop.

\item[{\rm (3)}] $\Ker(\tau)$ and $\tau(\mathbf{M})$ have the
disjunction property.
\end{enumerate}
\end{lemma}

\begin{proof}
(1) Suppose $\tau $ is faithful. Let  $F$ denote the least
nontrivial $\tau$-filter of $(\mathbf{M},\tau).$ There are two
cases: (i) If $\tau (M)\cap F\neq \left\{ 1\right\} $, then $\tau
(M)\cap F$ is the least nontrivial filter of $ \tau (\mathbf{M})$
and $\tau (\mathbf{M})$ is subdirectly irreducible. (ii) If $\tau
(\mathbf{M})\cap F=\left\{ 1\right\} $, then for all $x\in F$, $\tau
(x)=1$ because $\tau (x)\in \tau (M)\cap F$ and $F\subseteq
\Ker(\tau)=\left\{ 1\right\} $ is the trivial filter, a
contradiction. Therefore, only the first case is possible and
$\tau(\mathbf{M})$ is subdirectly irreducible.

Conversely, let $\tau(\mathbf{M})$ be subdirectly irreducible and
let $G$ be the least nontrivial filter of $\tau(\mathbf{M}).$ Then
the $\tau$-filter $F$ of $(\mathbf{M},\tau)$ generated by $G$ is the
least nontrivial $\tau$-filter of $(\mathbf{M},\tau).$ Indeed, if
$K$ is another nontrivial $\tau$-filter of $(\mathbf{M},\tau),$ then
$K\cap \tau(M) \supseteq F\cap \tau(M) = G.$  Then $K$ contains the
$\tau$-filter generated by $G$, that is $F\subseteq K$ which proves
$F$ is the least and $(\mathbf{M},\tau)$ is subdirectly irreducible.
\vspace{2mm}

Now let $(\mathbf{M},\tau)$ be subdirectly irreducible and let  $F$
denote the least nontrivial filter of $(\mathbf{M},\tau).$

(2) Suppose that $\tau$ is not faithful. Then $\Ker(\tau)$ is a
nontrivial $\tau $-filter. If $(\mathbf{M},\tau )$ is subdirectly
irreducible, it has a least nontrivial $\tau $-filter, $F$ say. So,
$F \subseteq \Ker(\tau)$, and hence $F$ is the least nontrivial
filter of the hoop $\Ker(\tau)$. Hence, $\Ker(\tau)$ is a
subdirectly irreducible hoop.

(3) Suppose, by way of contradiction, that for some $x\in
\Ker(\tau)$ and $y=\tau (y)\in \tau (M)$ one has $x<1$, $y<1$ and
$x\vee y=1$. It is easy to see that the BL-filters generated by $x$
and by $y$, respectively, are $\tau $-filters. Therefore they both
contain $F$. Hence, the intersection of these filters contains $F$.
Now let $c<1$ be in $F$. Then there is a natural number $n$ such
that $x^{n}\leq c$ and $y^{n}\leq c$. It follows that $ 1=(x\vee
y)^{n}=x^{n}\vee y^{n}\leq c$, a contradiction.
\end{proof}

\begin{lemma}\label{le:2.2}
If $(\mathbf{M},\tau )$ is a subdirectly irreducible state
BL-algebra, then $\tau (M)$ and $\Ker(\tau)$ are linearly ordered.
\end{lemma}

\begin{proof}
According to \cite[Thm 5.5]{CDH}, $\tau(M)$ is always linearly
ordered. On the other hand,  by Lemma \ref{le:2.1}, $\Ker(\tau)$ is
either a trivial hoop or a subdirectly irreducible  hoop, and hence
it is linearly ordered.
\end{proof}

\begin{theorem}\label{th:2.3}
Let $(\mathbf{M},\tau )$ be a state BL-algebra satisfying conditions
{\rm (1), (2)} and {\rm (3)} in Lemma {\rm \ref{le:2.1}}. Then $
(\mathbf{M},\tau )$ is subdirectly irreducible.
\end{theorem}

\begin{proof}
Suppose first that $\tau$ is faithful and that $\tau (\mathbf{M})$
is subdirectly irreducible. Let $F_{0}$ be the least nontrivial
filter of $\tau (\mathbf{M})$ and let $F$ be the BL-filter of
$\mathbf{M}$ generated by $F_{0}$. Then $F$ is a $\tau$-filter.
Indeed, if $x\in F$, then there is $\tau (a)\in F_{0}$ and a natural
number $n$ such that $\tau (a)^{n}\leq x$. It follows that $ \tau
(x)\geq \tau (\tau (a)^{n})=\tau (a)^{n}$, and $\tau (x)\in F$.

We assert that $F$ is the least nontrivial $\tau$-filter of
$(\mathbf{M},\tau )$. First of all, $\tau (\mathbf{M})$, being a
subdirectly irreducible BL-algebra, is linearly ordered. Now in
order to prove that $F$ is the least nontrivial $ \tau $-filter of
$(\mathbf{M},\tau )$, it suffices to prove that every $\tau $-filter
$G$ not containing $F$ is trivial. Now let $c<1$ in $F\backslash G$.
Then since $\Ker(\tau)=\left\{ 1\right\} $, $\tau (c)<1$. Next, let
$d\in G$. Then $ \tau (d)\in G$, and for every $n$ it cannot be
$\tau (d)^{n}\leq \tau (c)$, otherwise $\tau (c)\in G$. Hence, for
every $n$, $\tau (c)<\tau (d)^{n}$, and hence $\tau (c)$ does not
belong to the $\tau $-filter of $\tau (\mathbf{M})$ generated by
$\tau (d)$. By the minimality of $F$ in $\tau (\mathbf{M}),$ $\tau
(d)=1$ and since $\tau$ is faithful, we conclude that $d=1$ and $G $
is trivial, as desired.

Now suppose that $\Ker(\tau)$ is nontrivial. By condition (2),
$\Ker(\tau)$ is subdirectly irreducible. Thus, let $F$ be the least
nontrivial filter of $\Ker(\tau)$. Then $F$ is a non trivial $\tau
$-filter, and we have to prove that $F$ is the least nontrivial
$\tau $-filter of $(\mathbf{M},\tau )$. Let $G$ be any non trivial
$\tau $-filter of $(\mathbf{M},\tau )$. If $G\subseteq \Ker(\tau)$,
then it contains the least filter, $F$, of $\Ker(\tau)$, and
$F\subseteq G$. Otherwise, $G$ contains some $x\notin \Ker(\tau)$,
and hence it contains $\tau (x)<1$. Now by the disjunction property,
for all $y<1$ in $\Ker(\tau)$, $\tau (x)\vee y<1$ and $\tau (x)\vee
y\in \Ker(\tau)\cap G$. Thus, $G$ contains the filter generated by
$\tau (x)\vee y$, which is a non trivial filter of the hoop
$\Ker(\tau)$, and hence it contains $F$, the least nontrivial filter
of $ \Ker(\tau)$. This proves the claim.
\end{proof}

By \cite[Thm 3.5]{DKM},  conditions (1), (2), and (3) from Lemma
\ref{le:2.1} are independent ones even for state BL-algebras. In
addition, Theorem \ref{th:2.3} gives a characterization of
subdirectly irreducible state BL-algebras. If $(\mathbf{M},\tau)$ is
a state-morphism BL-algebra, combining \cite[Thm 4.5]{Dvu1} we can
say more about subdirectly irreducible state-morphism BL-algebras.
The following examples are from \cite{Dvu1}.

\begin{example}\label{ex:2}  {\rm
Let $\mathbf{M}$ be a  BL-algebra. On $M\times M$ we define two
operators, $\tau_1$ and $\tau_2,$ as follows
$$
\tau_1(a,b)=(a,a),\quad \tau_2(a,b)=(b,b),\quad (a,b)\in M\times
M.\eqno(2.0)
$$
Then $\tau_1$ and $\tau_2$ are two state-morphism operators on
$\mathbf M\times \mathbf M.$  Moreover, $(\mathbf{M}\times
\mathbf{M},\tau_1)$ and $(\mathbf{M}\times \mathbf{M},\tau_2)$ are
isomorphic state BL-algebras under the isomorphism $(a,b)\mapsto
(b,a).$}
\end{example}

We say that an element $a \in M$ is {\it Boolean} if $a^{--}=a$ and
$a\odot a =a.$  Let $B(\mathbf{M})$ be the set of Boolean elements.
Then $0,1\in B(\mathbf{M}),$ $B(\mathbf{M})$ is a subset of the
MV-skeleton $\MV(\mathbf{M}):= \{x \in M: x^{--}=x\},$ and $a\in
B(\mathbf{M})$ implies $a^-\in B(\mathbf{M}).$ We recall that
according to \cite[Thm 2]{TuSe}, $\MV(\mathbf{M})$ is an MV-algebra,
therefore, $B(\mathbf{M})$ is a Boolean subalgebra of
$\MV(\mathbf{M}).$

\begin{example}\label{ex:5.02} {\rm Let $\mathbf B$ be a local MV-algebra
such that $\Rad_1(\mathbf B)\ne \{1\}$ is a unique nontrivial filter
of $B.$ Let $\mathbf{M}$ be a subalgebra of $\mathbf B\times \mathbf
B$ that is generated by $\Rad_1(\mathbf B)\times \Rad_1(\mathbf B),$
that is $M=(\Rad_1(\mathbf B)\times \Rad_1(\mathbf B))\cup
(\Rad_1(\mathbf B)\times \Rad_1(\mathbf B))^-.$ Let
$\tau(x,y):=(x,x)$ for all $x,y \in M.$  Then $\tau$ is a
state-morphism operator on $\mathbf{M},$ $\Ker(\tau)=\{1\}\times
\Rad_1(\mathbf B) \subset \Rad_1(\mathbf{M})= \Rad_1(\mathbf
B)\times \Rad_1(\mathbf B),$ $\mathbf{M}$ has no Boolean nontrivial
elements, and $(\mathbf{M},\tau)$ is a subdirectly irreducible
state-morphism MV-algebra that is not linear.}
\end{example}

\begin{example}\label{ex:5.000} {\rm  Let $\mathbf A$ be a linear nontrivial
BL-algebra and $\mathbf B$ a nontrivial subdirectly irreducible
BL-algebra with the smallest nontrivial BL-filter $F_B$ and let $h:\
\mathbf A\to \mathbf B$ be a  BL-homomorphism.  On $\mathbf
\mathbf{M}=\mathbf A\times \mathbf B$ we define a mapping $\tau_h:\
M \to M$ by
$$ \tau_h(a,b) = (a,h(a)),\quad (a,b) \in M.\eqno(2.2)
$$
If we set $y=(0,1)$ and $y^-=(1,0)$, then  $y$ and $y^-$ are unique
nontrivial Boolean elements.

Then $\tau_h$ is a state-morphism operator on $\mathbf{M}$ and
$(\mathbf{M},\tau_h)$ is a subdirectly irreducible state-morphism
BL-algebra iff $\Ker(h)=\{a\in A: h(a)=1\}=\{1\}.$ In  such a case,
$\mbox{\rm Ker}(\tau_h) =\{1\}\times B$ and $F := \{1\} \times F_B$
is the least nontrivial state-morphism filter on $\mathbf{M}.$ }
\end{example}

Now we present the main result on the complete characterization of
subdirectly irreducible state-morphism BL-algebras which is a
combination of  \cite[Thm 4.5]{Dvu1} and Theorem \ref{th:2.3}.

\begin{theorem}\label{th:5.1}
A state-morphism BL-algebra $(\mathbf{M},\tau)$ is subdirectly
irreducible  if and only if one of the following three possibilities
holds.

\begin{itemize}

\item[{\rm (i)}]
$\mathbf M$ is linear,  $\tau =\Id_{M}$ is the identity on $M,$ and
the BL-reduct $\mathbf{M}$ is a subdirectly irreducible BL-algebra.

\item[{\rm (ii)}] The  state-morphism operator $\tau$ is not faithful,
$\mathbf{M}$ has no nontrivial Boolean elements, and the BL-reduct
$\mathbf{M}$ of $(\mathbf{M},\tau)$ is a local BL-algebra,
$\Ker(\tau)$ is a subdirectly irreducible irreducible hoop, and
$\Ker(\tau)$ and $\tau(\mathbf{M})$ have the disjunction property.

Moreover, $\mathbf{M}$ is linearly ordered if and only if
$\Rad_1(\mathbf{M})$ is linearly ordered, and in such a case,
$\mathbf{M}$ is a subdirectly irreducible BL-algebra such that if
$F$ is the smallest nontrivial state-filter for $(\mathbf{M},\tau),$
then $F$ is the smallest nontrivial BL-filter for $\mathbf{M}.$

If $\Rad_1(\mathbf{M})=\Ker(\tau),$ then $\mathbf{M}$ is linearly
ordered.

\item[{\rm (iii)}] The  state-morphism operator $\tau$ is not
faithful, $\mathbf{M}$ has a nontrivial  Boolean element. There are
a linearly ordered BL-algebra $\mathbf A,$ a subdirectly irreducible
BL-algebra $\mathbf B,$ and an injective BL-homomorphism $h:\
\mathbf A \to \mathbf B$ such that $(\mathbf{M},\tau)$ is isomorphic
as a state-morphism BL-algebra with the state-morphism BL-algebra
$(\mathbf A\times \mathbf B,\tau_h),$ where $\tau_h(x,y) =(x,h(x))$
for any $(x,y) \in A\times B.$
\end{itemize}
\end{theorem}

\begin{proof} It follows from \cite[Thm 4.5]{Dvu1} and Theorem
\ref{th:2.3}.
\end{proof}

We recall that a {\it t-norm} is a function $t:[0,1]\times [0,1]\to
[0,1]$ such that (i) $t$ is commutative, associative, (ii)
$t(x,1)=x,$ $x \in [0,1],$ and (iii) $t$ is nondecreasing in both
components. If $t$ is continuous, we define $x\odot_t y=t(x,y)$ and
$x\to_t y = \sup\{z\in [0,1]: t(z,x)\le y\}$ for $x,y \in [0,1],$
then $\mathbb I_t:=([0,1];\min,\max,\odot_t,\to_t,0,1)$ is a
BL-algebra. Moreover, according to \cite[Thm 5.2]{CEGT}, the variety
of all BL-algebras is generated by all $\mathbb I_t$ with a
continuous t-norm $t.$ Let $\mathcal T$ denote the system of all
BL-algebras $\mathbb I_t,$ where $t$ is any continuous t-norm.

The proof of the following result will follow from Theorem
\ref{ad:5.2}.

\begin{theorem}\label{th:3.4}
The variety of all state-morphism BL-algebras is generated by the
system $\{D(\mathbb I_t): t \in \mathcal T\}.$
\end{theorem}

\section{General State-Morphism Algebras}%3

In this section, we generalize the notion of state-morphism
BL-algebras to an arbitrary variety of algebras of some type. It is
interesting that many results known only for state-morphism
MV-algebras or state-morphism BL-algebras have a very general
presentation as state-morphism algebras. The main result of this
section, Theorem \ref{diag}, says that every subdirectly irreducible
state-morphism algebra can be embedded into some diagonal one.

Let   $\mathbf A$ be any algebra of type $F$ and let
$\mathrm{Con}\,\mathbf A$ be the system of congruences on $\mathbf
A$ with the least congruence $\Delta_{\mathbf A}.$ An endomorphisms
$\tau:\mathbf A\longrightarrow \mathbf A$ satisfying $\tau\circ\tau
= \tau$ is said to be a {\it state-morphism} on $\mathbf A$ and a
couple $(\mathbf A,\tau)$ is said to be a {\it state-morphism
algebra} or an algebra with internal state-morphism. Clearly, if
$\mathcal K$ is a variety of algebras of type $F$, then the class
$\mathcal K_\tau$ of all state-morphism algebras $(\mathbf A,
\tau)$, where $\mathbf A\in\mathcal K$ and $\tau$ is any
state-morphism on $\mathbf A,$ forms a variety, too.

In the rest of the paper,   we will assume that $\mathbf A$  is an
arbitrary algebra with a fixed type $F;$ if $\mathbf A$ is of a
specific type, it will be said that and specified.

\begin{definition}{\rm
Let $\mathbf B\in\mathcal K$. Then an algebra $D(\mathbf
D):=(\mathbf B\times \mathbf B,\tau_B),$ where $\tau_B$ is defined
by $\tau_B(x,y)=(x,x),$ $x,y \in B,$ is a state-morphism algebra
(more precisely $(\mathbf B\times \mathbf B,\tau_B)\in\mathcal
K_\tau$); we call $\tau_B$ also a {\it diagonal state-operator}. If
a state-morphism algebra $(\mathbf C,\tau)$ can be embedded into
some diagonal state-morphism algebra, $(\mathbf{B}\times
\mathbf{B},\tau_B),$  $(\mathbf C,\tau)$ is said to be a {\it
subdiagonal} state-morphism algebra}, or, more precisely, $\mathbf
B$-subdiagonal.
\end{definition}

Let  $(\mathbf A,\tau)$ be a state-morphism algebra. We introduce
the following sets:
 $$\theta_\tau =\{(x,y)\in A\times A: \tau (x)=\tau(y)\},\eqno(3.1)$$
 $$\tau (A) = \{\tau (x): x\in A\}.
 $$
The subalgebra which is an image of $\mathbf A$ by $\tau$ is denoted
by $\tau (\mathbf A)$ and thus $\tau (\mathbf A)\in\mathcal K$ and
$(\tau(\mathbf A),\Id_{\tau(A)}) \in \mathcal K_\tau,$ where
$\Id_{\tau(A)}$ is the identity on $\tau(A);$ we have also
$\tau|\tau(A)=\Id_{\tau(A)}.$

If $\phi\in\mathrm{Con}\,\tau(\mathbf A),$ we define
$$
\theta_\phi :=\{(x,y)\in A\times A: (\tau (x),\tau
(y))\in\phi\}.\eqno(3.2)
$$

Finally, if $\phi\subseteq A^2$ then the congruence on $\mathbf A$
generated by $\phi$ is denoted by $\Theta (\phi)$ and the congruence
on $(\mathbf A,\tau)$  generated by $\phi$ is denoted by
$\Theta_\tau (\phi)$. Clearly $\mathrm{Con}\, (\mathbf A,\tau)
\subseteq \mathrm{Con}\,\mathbf A$ and also $\Theta (\phi)\subseteq
\Theta_\tau (\phi)$.

\begin{lemma}\label{1}
Let  $(\mathbf A,\tau)$ be a state-morphism algebra. For any
$\phi\in\mathrm {Con}\,\tau (\mathbf A),$ we have
$\theta_\phi\in\mathrm{Con}\, (\mathbf A,\tau), $ and $\theta_\phi
\cap \tau(A)^2 = \phi.$ In addition, $\theta_\tau\in\mathrm{Con}\,
(\mathbf A,\tau),$  $\phi\subseteq \theta_\phi,$ and
$\Theta_\tau(\phi) \subseteq \theta_\phi.$
\end{lemma}

\begin{proof}
Clearly, $\theta_\phi$ is reflexive and symmetric. Moreover, if
$(x,y),(y,z)\in\theta_\phi$, then
$\bigl(\tau(x),\tau(y)\bigr),\bigl(\tau(y),\tau (z)\bigr)\in\phi$
and thus $\bigl(\tau (x),\tau(z)\bigr)\in \phi$ which gives
$(x,z)\in\theta_\phi$.

Let $f^{\mathbf A}$ be an $n$-ary operation on $\mathbf A$ and let
$x_1,\ldots,x_n,y_1,\ldots, y_n\in A$ be such that
$(x_i,y_i)\in\theta_\phi$ for any $i=1,\ldots, n$. Then $\bigl(\tau
(x_i),\tau (y_i)\bigr)\in\phi$ holds for any $i=1,\ldots,n$. Due to
$\phi\in\mathrm{Con}\,\tau (\mathbf A),$ we obtain $\bigl(f^{\tau
(\mathbf A)}(\tau (x_1),\ldots, \tau (x_n)),f^{\tau (\mathbf
A)}(\tau (y_1),\ldots, \tau (y_n))\bigr)\in\phi$.

Because $\tau$ is an endomorphism,   $\tau (f^\mathbf A (x_1,\ldots,
x_n)) = f^{\tau (\mathbf A)}(\tau (x_1),\ldots, \tau (x_n))$ and
$\tau (f^\mathbf A (y_1,\ldots, y_n)) = f^{\tau (\mathbf A)}(\tau
(y_1),\ldots, \tau (y_n))$ which gives $\bigl(\tau (f^\mathbf A
(x_1,\ldots, x_n)),$ $\tau (f^\mathbf A (y_1,\ldots,
y_n))\bigr)\in\phi$ and finally also $ \bigl(f^\mathbf A
(x_1,\ldots, x_n),f^\mathbf A (y_1,\ldots,
y_n)\bigr)\in\theta_\phi$.

Moreover, take  an arbitrary $(x,y)\in\theta_\phi.$ Then
$\bigl(\tau(\tau (x)),\tau(\tau (y))\bigr)=\bigl(\tau (x),\tau
(y)\bigr)\in\phi$ which gives $\bigl(\tau (x),\tau
(y)\bigr)\in\theta_\phi$.

Hence, $\theta_\phi\in\mathrm{Con}\,(\mathbf A,\tau)$ and if
$\phi=\Delta_{\tau (\mathbf A)},$ then $\theta_\phi=\theta_\tau$.

It is clear that $\theta_\phi \cap \tau(A)^2 \supseteq \phi.$  Now
let $(x,y) \in \theta_\phi \cap \tau(A)^2.$ Then $x,y\in \tau(A),$
$(\tau(x),\tau(y)) \in \phi \subseteq \tau(A)^2,$ so that $x=\tau(x)
\in \tau(A),$ $y=\tau(y) \in \tau(A),$ and consequently, $(x,y) \in
\phi.$

It is evident that $\theta_\tau$ is a congruence on $(\mathbf
A,\tau).$

Finally, if $(x,y)\in\phi$ then $\tau (x)=x$ and $\tau (y)=y$ which
gives $\bigl(\tau (x),\tau (y)\bigr)=(x,y)\in\phi$. Thus
$(x,y)\in\theta_\phi$ which finishes the proof that $\phi\subseteq
\theta_\phi$ and $\Theta_\tau(\phi) \subseteq \theta_\phi.$
\end{proof}

\begin{lemma}\label{2}
Let $\theta\in\mathrm{Con}\,\mathbf A$ be such that
$\theta\subseteq\theta_\tau$. Then $\theta\in\mathrm{Con}\,(\mathbf
A,\tau)$ holds.

Moreover, if $x,y\in A$ are such that $(x,y)\in\theta_\tau,$ then
$\Theta (x,y)=\Theta_\tau (x,y)$.
\end{lemma}

\begin{proof}
If $(x,y)\in\theta\subseteq\theta_\tau,$  then $\tau (x)=\tau (y)$
and thus $\bigl(\tau (x),\tau (y)\bigr)=\bigl(\tau (x),\tau
(x)\bigr)\in\theta$ proves that $\theta\in\mathrm {Con}\,(\mathbf
A,\tau)$.

Moreover, if $(x,y)\in\theta_\tau,$ then $\Theta (x,y)\subseteq
\theta_\tau$. Due to the first part of Lemma, we obtain $\Theta
(x,y)\in\mathrm{Con}\,(\mathbf A,\tau)$ and thus $\Theta_\tau
(x,y)\subseteq\Theta (x,y)$ holds. The second inclusion is trivial.
\end{proof}

\begin{lemma}\label{3}
If $x,y\in\tau (\mathbf A),$ then $\Theta (x,y)=\Theta_\tau (x,y).$ Consequently,
$\Theta(\phi) = \Theta_\tau(\phi)$ whenever $\phi \subseteq \tau(A)^2.$

\end{lemma}

\begin{proof}
Let us denote by $\phi$ the congruence on $\tau (\mathbf A)$
generated by $(x,y)$. Clearly we obtain the chain of inclusions
$\phi\subseteq\Theta (x,y)\subseteq\Theta
(\phi)\subseteq\theta_\phi$ (because $(x,y)\in\phi$ and
$\phi\subseteq \theta_\phi$, see Lemma \ref{1}).

Assume $(a,b)\in\Theta (x,y),$ then $(a,b)\in\theta_\phi$ and
thus $(\tau (a),\tau (b))\in\phi\subseteq\Theta (x,y)$. Thus $\Theta
(x,y)\in\mathrm {Con}\,(\mathbf A,\tau)$ and $\Theta_\tau
(x,y)\subseteq\Theta (x,y)$ holds. The second inclusion is trivial.

Finally, let $\phi \subseteq \tau(A)^2.$ By \cite[Thm 5.3]{BuSa},
the both congruence lattices of $\mathbf A$ and $(\mathbf A,\tau)$
are complete sublattices of the lattice of equivalencies on $\mathbf
A,$ and therefore, they have the same infinite suprema. Hence, by
the first part of the lemma,
$$
\Theta(\phi) = \bigvee_{(x,y)\in\phi} \Theta(x,y) =
\bigvee_{(x,y)\in\phi} \Theta_\tau(x,y)=\Theta_\tau(\phi).
$$
\end{proof}

\begin{remark}\label{re:1} {\rm By Lemma \ref{1},
if $\phi$ is a congruence on $\tau(\mathbf A),$ then $\theta_\phi$
is an extension of $\phi$ on $(\mathbf A,\tau)$ and $\Theta(\phi)
=\Theta_\tau(\phi) \subseteq \theta_\phi.$ There is a natural
question whether  $\Theta(\phi) = \theta_\phi$ ?  The answer is
positive if and only if $\tau$ is the identity on $A$.   Indeed, if
$\tau$ is the identity on $A,$ the statement is evident, in the
opposite case, we have $\theta_{\Delta_{\tau(\mathbf
A)}}=\theta_\tau\not=\Delta_\mathbf A = \Theta(\Delta_{\tau(\mathbf
A)}).$ }
\end{remark}

\begin{theorem}\label{4}
Let $(\mathbf A,\tau)$ be a subdirectly irreducible state-morphism
algebra such that $\mathbf A$ is  subdirectly reducible. Then there
is a subdirectly irreducible algebra $\mathbf B$ such that $(\mathbf
A,\tau)$ is $\mathbf B$-subdiagonal.
\end{theorem}

\begin{proof}
First, if $\theta_\tau = \Delta_{\mathbf A},$ then for any $x\in A,$
the equality $\tau (x)= x$ holds and thus $\mathrm {Con}\,\mathbf A
= \mathrm {Con}\,(\mathbf A,\tau)$ which is absurd because $\mathbf
A$ is subdirectly irreducible and $(\mathbf A,\tau)$ is not
subdirectly irreducible.

The subdirect irreducibility of $(\mathbf A,\tau )$ implies that
there is a least proper congruence $\theta_{\min}\in\mathrm
{Con}\,(\mathbf A,\tau)$. Moreover, due to Lemma \ref{2}, the
congruence $\theta_{\min}$ is also a least proper congruence
$\theta$ on $\mathbf A$ with $\theta\subseteq\theta_\tau$  and thus
$\theta_{\min}$ is an atom in $\mathrm {Con}\,\mathbf A$. Let us
denote
$$
\theta_\tau^{\perp}=\{\theta\in\mathrm {Con}\,\mathbf A:
\theta\cap\theta_\tau=\Delta_{\mathbf A}\}.
$$
First, we prove that there exists proper
$\theta\in\theta_{\tau}^{\perp}$. The subdirect reducibility of
$\mathbf A$ shows that there exists proper $\theta\in\mathrm
{Con}\,\mathbf A$ with $\theta_{\min}\not\subseteq\theta$. Hence,
$\theta_\tau\cap\theta=\Delta_{\mathbf A}$ holds (because if
$\theta_\tau\cap\theta\not=\Delta_{\mathbf A},$ then
$\theta_\tau\cap\theta$ is a proper congruence contained in
$\theta_\tau$ and minimality of $\theta_{\min}$ yields
$\theta_{\min}\subseteq\theta\cap\theta_\tau\subseteq\theta,$ which
is absurd).

Moreover, let us have $\theta_n\in\theta_\tau^{\perp}$  for any
$n\in\mathbb N$ with  $\theta_n\subseteq\theta_{n+1},$ then clearly
$\bigvee_{n\in\mathbb N}\theta_n=\bigcup_{n\in\mathbb
N}\theta_n\in\theta_\tau^{\perp}$. Due to Zorn's Lemma, there is
maximal $\theta^*\in\theta_\tau^{\perp}$.

We have proved that both $\theta_\tau$ and $\theta^*$ are proper
congruences on $\mathbf A$ with
$\theta_\tau\cap\theta^*=\Delta_{\mathbf A}.$ By the Birkhoff
Theorem about subdirect reducibility, $\mathbf A$ is a subdirect
product of two algebras $\mathbf A /\theta_\tau$ and $\mathbf
A/\theta^*$ with an embedding $h:\mathbf A\longrightarrow \mathbf
A/\theta_\tau\times\mathbf A/\theta^*$ defined by
$h(x)=(x/\theta_\tau,x/\theta^*).$

Now we define the mapping $\psi: A/\theta_\tau\longrightarrow
A/\theta^*$ by $\psi (x/\theta_\tau)=\tau (x)/\theta^*.$ Clearly
$\psi$ is well-defined because $x/\theta_\tau=y/\theta_\tau$ yields
$\tau (x)=\tau (y)$ and thus $\psi (x/\theta_\tau)=\tau (x)/\theta^*
=\tau (y)/\theta^*=\psi (y/\theta_\tau).$ Let us suppose that $\psi
(x/\theta_\tau)=\psi (y/\theta_\tau).$ Then $\tau (x)/\theta^*=\tau
(y)/\theta^*$ and $\bigl(\tau (x),\tau (y)\bigr)\in\theta^*$. Hence,
$\Theta (\tau (x),\tau (y))\subseteq \theta^*$ holds. Finally, if
$\tau (x)\not =\tau (y)$ (thus $\Theta (\tau (x),\tau (y))$ is a
proper congruence), then $\tau (x),\tau (y)\in\tau (\mathbf A)$ and
Lemma \ref{3} yields $\Theta (\tau (x),\tau (y))\in\mathrm {Con}\,
(\mathbf A,\tau)$ and thus $\theta_{\min}\subseteq\Theta (\tau
(x),\tau (y))\subseteq\theta^*$ which is absurd
($\theta_{\min}\subseteq \theta_\tau\cap\theta^*=\Delta_{\mathbf
A}$). Therefore, the mapping $\psi$ is injective.

We shall prove that $\psi$ is a homomorphism (and thus an
embedding). If $f^{\mathbf A}$ is an $n$-ary operation and
$x_1/\theta_\tau,\ldots, x_n/\theta_\tau\in\mathbf A/\theta_\tau,$
then
\begin{eqnarray*}
\psi (f^{\mathbf A/\theta_\tau}(x_1/\theta_\tau ,\ldots , x_n/\theta_\tau))&=& \psi (f^\mathbf A (x_1,\ldots ,x_n)/\theta_\tau)\\
&=& \tau (f^{\mathbf A}(x_1,\ldots ,x_n))/\theta^*\\
&=& f^{\mathbf A}(\tau (x_1),\ldots , \tau (x_n))/\theta^*\\
&=& f^{\mathbf A/\theta^*}(\tau (x_1)/\theta^* ,\ldots, \tau (x_n)/\theta^*)\\
&=& f^{\mathbf A/\theta^*}(\psi (x_1/\theta_\tau ),\ldots ,\psi
(x_n/\theta_\tau )).
\end{eqnarray*}

Now we prove that $\mathbf A$ is $\mathbf A/\theta^*$-diagonal. Let
$g: A\longrightarrow (A/\theta^*)^2$ be defined via $g(x)=(\psi(
x/\theta_\tau),x/\theta^*)=(\tau (x)/\theta^*,x/\theta^*)$. Because
the mapping $g$ is the composition of two functions $h$ and $\psi$
which are embeddings,  $g$ is also an embedding of $\mathbf A$ into
$(\mathbf A/\theta^*)^2$. Now we can compute:
\begin{eqnarray*}
g(\tau (x))&=& \bigl(\tau (\tau (x))/\theta^*,\tau (x)/\theta^*\bigr)\\
&=& \bigl(\tau (x)/\theta^*,\tau (x)/\theta^*\bigr)\\
&=& \tau_{\mathbf
A/\theta^*} (\tau (x)/\theta^* , x/\theta^*)\\
&=& \tau_{\mathbf A/\theta^*} (g(x)),
\end{eqnarray*}
where $\tau_{\mathbf A/\theta^*}$ is the diagonal state-morphism on
the product $\mathbf A/\theta^* \times \mathbf A/\theta^*.$
Therefore, $g:(\mathbf A,\tau)\longrightarrow (\mathbf
A/\theta^*\times \mathbf A/\theta^*,\tau_{\mathbf A/\theta^*})$ is
an embedding and $(\mathbf A,\tau)$ is  $\mathbf
A/\theta^*$-diagonal.

Finally, we prove the subdirect irreducibility of $\mathbf
A/\theta^*$. Of course, $\theta_{\min}\cap\theta^*=\Delta_{\mathbf
A}$ yields $\theta_{\min}\not\subseteq\theta^*$ and thus
$\theta^*\subset \theta^*\vee \theta_{\min}$. Moreover, if
$\theta^*\subset \theta,$ from maximality of $\theta^*$ we obtain
$\theta\cap\theta_\tau\not=\Delta_{\mathbf A}$ and thus
$\theta_{\min}\subseteq \theta_\tau\cap\theta$. Finally,
$\theta_{\min}\vee\theta^*\subseteq
(\theta_\tau\cap\theta)\vee\theta^*\subseteq
(\theta_\tau\cap\theta)\vee\theta =\theta$ holds. Hence, for any
congruence $\theta\in\mathrm {Con}\,\mathbf A,$ the inequality
$\theta^*\subset \theta^*\cap\theta_{\min}\subseteq\theta$ holds.
Due to the Birkhoff's Theorem and the Second Homomorphism Theorem,
an algebra $\mathbf A/\theta^*$ is subdirectly irreducible.
\end{proof}

Theorem \ref{4}  can be extended as follows.

\begin{theorem}\label{diag}
For every subdirectly irreducible state-morphism algebra $(\mathbf
A,\tau),$  there is a subdirectly irreducible algebra $\mathbf B$
such that $(\mathbf A,\tau)$ is $\mathbf B$-subdiagonal.
\end{theorem}

\begin{proof}
There are two cases: (1)  $(\mathbf A,\tau)$ and $\mathbf A$ are
subdirectly irreducible, and (2) $(\mathbf A,\tau)$ is a subdirectly
irreducible state-morphism algebra and $\mathbf A$ is a subdirectly
reducible algebra

(1) Assume that $(\mathbf A,\tau)$ and $\mathbf A$ are subdirectly
irreducible. Define two state-morphism algebras $(\tau(\mathbf
A)\times \mathbf A,\tau_1)$ and $(\mathbf A \times \mathbf
A,\tau_2),$ where $\tau_1(a,b)= (a,a),$ $(a,b) \in \tau(A)\times A,$
and $\tau_2(a,b) =(a,a),$ $a,b \in A.$ Then the first one is a
subalgebra of the second one.

Define a mapping $\phi: A \to \tau(A)\times A$ defined by
$\phi(a)=(\tau(a),a),$ $a\in A.$  Then $\phi$ is injective because
if $\phi(a)=\phi(b)$ then $(\tau(a),a)=(\tau(b),b)$ and $a=b.$ We
show that $\phi$ is a homomorphism. Let $f^{\mathbf A}$ be an
$n$-ary operation on $\mathbf A$ and let $a_1,\ldots, a_n \in A.$
Then
\begin{eqnarray*} \phi(f^{\mathbf A}(a_1,\ldots,a_n)) &=& \bigl(\tau(f^{\mathbf
A}(a_1,\ldots,a_n)),f^{\mathbf A}(a_1,\ldots,a_n)\bigr)\\
&=& \bigl(f^{\mathbf A}(\tau(a_1),\ldots,\tau(a_n)),f^{\mathbf
A}(a_1,\ldots,a_n) \bigr)\\
&=& f^{\tau(\mathbf A)\times \mathbf A}\bigl((\tau(a_1),a_1),\ldots, (\tau(a_n),a_n)\bigr)\\
&=& f^{\tau(\mathbf A)\times \mathbf A}(\phi(a_1),\ldots,\phi(a_n)).
\end{eqnarray*}

Since $\phi: \mathbf A \to \tau(\mathbf A) \times \mathbf A
\subseteq \mathbf A \times \mathbf A,$ $\phi$ can be assumed also as
an injective homomorphism from the state-morphism algebra $(\mathbf
A,\tau)$ into the state-morphism algebra $ D(\mathbf B),$ where
$\mathbf B:= \mathbf A$ is a subdirectly irreducible algebra.

(2)  This case was proved in Theorem \ref{4}.
\end{proof}

For example, a state-morphism algebra $(\mathbf A,\Id_A),$ where
$\Id_A$ is the identity on $A,$ is subdirectly irreducible if and
only if $\mathbf A$ is subdirectly irreducible. Therefore, $(\mathbf
A,\Id_A)$ can be embedded into $(\mathbf A \times \mathbf A,
\tau_A)$ under the mapping $a \mapsto (a,a),$ $a \in A.$
Consequently, every subdirectly irreducible state-morphism algebra
$(\mathbf A,\Id_A)$ is $\mathbf A$-subdiagonal with $\mathbf A$
subdirectly irreducible.

We note that in the same way as in \cite[Lem 6.1]{DKM}, it is
possible to show that the class of subdiagonal state-morphism
algebras is closed under subalgebras and ultraproducts, and not
closed under homomorphic images, see \cite[Lem 6.6]{DKM}.

\section{Varieties of State-Morphism Algebras and Their
Generators}%4

In this section, we study varieties of state-morphism algebras and
their generators. It is interesting to note that some similar
results proved for state-morphism MV-algebras in \cite{DKM} can be
obtained in an analogous way also for a general variety of algebras.

Let $\tau$ be a state-morphism operator on an algebra $\mathbf A.$
We set
$$\Ker(\tau):=\{(x,y)\in A\times A: \tau(x)=\tau(y)\},
$$
the {\it kernel} of $\tau.$  We say that $\tau$ is {\it faithful} if
$\Ker(\tau)=\Delta_{\mathbf A}.$ It is evident that $\tau$ is
faithful iff $\tau(x)=x$ for each $x \in A.$ In addition, $\tau$ is
faithful iff $\tau$ is injective.

For every class $\mathcal{K}$ of same type algebras, we set
$\mathsf{D}(\mathcal{K})=\left\{ D(\mathbf{A}):\mathbf{A}\in
\mathcal{K} \right\},$ where
$D(\mathcal{\mathbf{A}})=(\mathbf{A}\times \mathbf{A},\tau _{A}).$

As usual, given a class $\mathcal{K}$ of algebras of the same type,
$\mathsf{ I}(\mathcal{K})$, $\mathsf{H}(\mathcal{K})$,
$\mathsf{S}(\mathcal{K})$ and $ \mathsf{P}(\mathcal{K})$ and
$\mathsf P _{\mathsf U}(\mathcal K)$ will denote the class of
isomorphic images, of homomorphic images, of subalgebras, of direct
products and of ultraproducts of algebras from $ \mathcal{K}$,
respectively. Moreover, $\mathsf{V}(\mathcal{K})$ will denote the
variety generated by $\mathcal{K}$.

\begin{lemma}\label{wd1}
{\rm (1)} Let $\mathcal{K}$ be a class of algebras of the same type
$F$. Then $\mathsf{VD}(\mathcal{K})\subseteq
\mathsf{V}(\mathcal{K})_\tau$.
\newline {\rm (2)} Let $\mathcal{V}$
be any variety. Then $\mathcal{V} _\tau=\mathsf{ISD}(\mathcal{V})$.
\end{lemma}

\begin{proof}
(1) If $ D(\mathbf{A})\in\mathsf D (\mathcal K)$ (thus
$\mathbf{A}\in\mathcal K$), then the $F$-reduct of the algebra
$D(\mathbf{A})$ is the algebra $\mathbf{A}\times \mathbf{A}$ which
belongs to the variety $\mathsf V(\mathcal K)$. Due to definition of
$\mathsf V (\mathcal K)_\tau,$ we obtain also
$D(\mathbf{A})\in\mathsf V (\mathcal K)_\tau$. We have proved that
$\mathsf D (\mathcal K)\subseteq \mathsf V (\mathcal K)_\tau$.
Because $\mathsf V (\mathcal K)_\tau$ is a variety then also
$\mathsf {VD} (\mathcal K)\subseteq \mathsf V (\mathcal K)_\tau$

(2) Let $(\mathbf{A},\tau )\in \mathcal{V}_\tau$. As we have seen in
the proof of Theorem \ref{diag}, the map $ \phi :a\mapsto (\tau
(a),a)$ is an injective homomorphism of $(\mathbf{A},\tau )$ into
$D( \mathbf{A})$. Hence, $\phi $ is compatible with $\tau $, and
$(\mathbf{A},\tau )\in \mathsf{ISD}(\mathcal{V})$. Conversely, the
$F$-reduct of any algebra in $\mathsf{D}(\mathcal{V})$ is in
$\mathcal{V}$, (being a direct product of algebras in
$\mathcal{V}$), and hence the $F$-reduct of any member of
$\mathsf{ISD}(\mathcal{V})$ is in $\mathsf{IS}(\mathcal{V})=
\mathcal{V}$. Hence, any member of $\mathsf{ISD}(\mathcal{V})$ is in
$\mathcal{V }_\tau$.
\end{proof}

\begin{lemma}\label{main1}
Let $\mathcal{K}$ be a class of algebras of the same type $F$. Then:
\newline {\rm (1)}
$\mathsf{DH}(\mathcal{K})\subseteq \mathsf{HD}(\mathcal{K})$.
\newline {\rm (2)}
$\mathsf{DS}(\mathcal{K})\subseteq \mathsf{ISD}(\mathcal{K})$.
\newline {\rm (3)}
$\mathsf{DP}(\mathcal{K})\subseteq \mathsf{IPD}(\mathcal{K})$.
\newline {\rm (4)}
$\mathsf{VD}(\mathcal{K})=\mathsf{ISD}(\mathsf{V}( \mathcal{K}))$.
\end{lemma}

\begin{proof}
(1) Let $D(\mathbf C)\in \mathsf{DH}(\mathcal{K})$. Then there are $
\mathbf{A}\in \mathcal{K}$ and a homomorphism $h$ from $\mathbf{A}$
onto $ \mathbf C$. Let for all $a,b\in A$, $h^{*}(a,b)=(h(a),h(b))$.
We claim that $h^{*}$ is a homomorphism from $D(\mathbf{A})$ onto
$D(\mathbf C)$. That $ h^{*}$ is a homomorphism is clear. We verify
that $h^{*}$ is compatible with $\tau _{A}$. We have $h^{*}(\tau
_{A}(a,b))=h^{*}(a,a)=(h(a),h(a))=\tau _{C}(h(a),h(b))=\tau
_{C}(h^{*}(a,b)).$ Finally, since $h$ is onto, given $(c,d)\in
C\times C$, there are $a,b\in A$ such that $h(a)=c$ and $h(b)=d$.
Hence, $h^{*}(a,b)=(c,d)$, $h^{*}$ is onto, and $D(\mathbf C)\in
\mathsf{HD}(\mathcal{K})$.

(2) It is trivial.

(3) Let $\mathbf{A}=\prod_{i\in I}(\mathbf{A}_{i})\in
\mathsf{P}(\mathcal{K})$, where each $\mathbf{A}_{i}$ is in
$\mathcal{K}$. Then the map
\begin{center}
$\Phi :\bigl((a_{i}:i\in I),(b_{i}:i\in I)\bigr)\mapsto
\bigl((a_{i},b_{i}\bigr):i\in I)$
\end{center}
is an isomorphism from $D(\mathbf{A})$ onto $\prod_{i\in
I}D(\mathbf{A}_{i})$. Indeed, it is clear that $\Phi $ is an
$F$-isomorphism. Moreover, denoting the state-morphism of
$\prod_{i\in I}D(\mathbf{A}_{i})$ by $\tau ^{*}$, we get:
\begin{eqnarray*}&\Phi \bigl(\tau _{A}\bigl((a_{i}:i\in I),(b_{i}:i\in I)\bigr)
\bigr)=\Phi
\bigl((a_{i}:i\in I),(a_{i}:i\in I)\bigr)=\\
&= \bigl((a_{i},a_{i}):i\in I\bigr)=\bigl(\tau
_{\mathbf{A}_{i}}(a_{i},b_{i}):i\in I\bigr)=\tau ^{*}\bigl(\Phi
\bigl((a_{i}:i\in I),(b_{i}:i\in I)\bigr)\bigr),
\end{eqnarray*}
and hence $\Phi $ is an isomorphism.

(4) By (1), (2) and (3),
$\mathsf{DV}(\mathcal{K})=\mathsf{DHSP}(\mathcal{K} )\subseteq
\mathsf{HSPD}(\mathcal{K})=\mathsf{VD}(\mathcal{K})$, and hence $
\mathsf{ISDV}(\mathcal{K})\subseteq
\mathsf{ISVD}(\mathcal{K})=\mathsf{VD}( \mathcal{K})$. Conversely,
by Lemma \ref{wd1}(1), $\mathsf{VD}(\mathcal{K} )\subseteq
\mathsf{V}(\mathcal{K})_\tau$, and by Lemma \ref{wd1}(2), $
\mathsf{V}(\mathcal{K})_\tau=\mathsf{ISDV}(\mathcal{K})$. This
proves the claim.
\end{proof}

\begin{theorem}\label{th:3.3.1}
{\rm (1)} For every class $\mathcal K$ of algebras of the same type
$F,$  $\mathsf{V}(\mathsf{D}(\mathcal K))=\mathsf{V}( \mathcal
K)_\tau$.

{\rm (2)} Let $\mathcal K_1$ and $\mathcal K_2$ be two classes of
same type algebras. Then $\mathsf{V}(D( \mathcal
K_1))=\mathsf{V}(D(\mathcal K_2))$ if and only if
$\mathsf{V}(\mathcal K_1)=\mathsf{V }(\mathcal K_2)$.

\end{theorem}

\begin{proof}
(1) By Lemma \ref{main1}(4), $\mathsf{VD}(\mathcal
K)=\mathsf{I}$\textsf{S}$\mathsf{D}(\mathsf{V}(\mathcal K))$.
Moreover, by Lemma \ref{wd1}(2), $\mathsf{V}(\mathcal
K)_\tau=\mathsf{ISDV }(\mathcal K)$. Hence,
$\mathsf{V}(\mathsf{D}(\mathcal K))=\mathsf{V}(\mathcal K )_\tau$.

(2) We have $\mathsf{V}(\mathsf{D}(\mathcal
K_1))=\mathsf{V}(\mathcal K_1)_\tau$ and $
\mathsf{V}(\mathsf{D}(\mathcal K_2))=\mathsf{V}(\mathcal K_2)_\tau$.
Clearly, $\mathsf{V} (\mathcal K_1)=\mathsf{V}(\mathcal K_2)$
implies $\mathsf{V}(\mathcal K_1)_\tau= \mathsf{V}(\mathcal
K_2)_\tau$, and hence $\mathsf{V}(\mathsf{D}(\mathcal
K_1))=\mathsf{V }(\mathsf{D}(\mathcal K_2))$. Conversely,
$\mathsf{V}(\mathsf{D}(\mathcal K_1))=\mathsf{V}(\mathsf{D}(\mathcal
K_2))$ implies $\mathsf{V}(\mathcal K_1)_\tau=\mathsf{V}(\mathcal
K_2)_\tau$. But any algebra $\mathbf{A}\in \mathsf{V}(\mathcal K_1)$
is the $F$-reduct of a state-morphism algebra in
$\mathsf{V}(\mathcal K_1)_\tau$, namely of $(\mathbf{A}, \Id_A).$

It follows that, if $\mathsf{V}(\mathcal K_1)_\tau=
\mathsf{V}(\mathcal K_2)_\tau$, then the classes of $F$-reducts of
$\mathsf{V}(\mathcal K_1)_\tau$ and of $\mathsf{V}(\mathcal K_2
)_\tau$ coincide, and hence $\mathsf{V}(\mathcal
K_1)=\mathsf{V}(\mathcal K_2)$.
\end{proof}

As a direct corollary of Theorem \ref{th:3.3.1}, we have:

\begin{theorem}\label{th:3.4.1}
If a system $\mathcal K$ of algebras of the same type $F$ generates
the whole variety $\mathcal V(F)$ of all algebras of type $F,$ then
the variety $\mathcal V(F)_\tau$ of all state-morphism algebras
$(\mathbf A,\tau),$ where $\mathbf  A \in \mathcal V(F),$ is
generated by the class $\{D(\mathbf A): \mathbf A\in \mathcal K\}.$
\end{theorem}

Some applications of the latter theorem for different varieties of
algebras will be done in Section 5.

\begin{theorem}
If $\mathbf A$ is a subdirectly irreducible algebra, then any
state-morphism algebra $(\mathbf A,\tau)$ is subdirectly
irreducible.
\end{theorem}

\begin{proof}
Let $\mathbf A$ be a subdirectly irreducible algebra and let $\tau$
be a state-morphism operator on $\mathbf A.$ If $\tau$ is the
identity on $A$,  then $\mathrm{Con}\, \mathbf A = \mathrm{Con}\,
(\mathbf A,\tau)$ and, consequently,  $(\mathbf A,\tau)$ is
subdirectly irreducible. If $\tau$ is not the identity on $A,$ then
$\theta_\tau,$ defined by (3.1),  is a nontrivial congruence on
$\mathbf A,$ and thus $\theta_{\min}\subseteq\theta_\tau,$ where
$\theta_{\min}\in\mathrm{Con}\,\mathbf A$ is the least nontrivial
congruence. Hence, $\theta_{\min}$ belongs  to the set
$\mathrm{Con}\,(\mathbf A,\tau)$, see Lemma \ref{2}. Therefore,
$\mathrm{Con}\,(\mathbf A,\tau)\subseteq \mathrm{Con}\,\mathbf A$
yields the subdirect irreducibility of the algebra $(\mathbf
A,\tau)$, more precisely, $\theta_{\min}$ is also the least proper
congruence in $\mathrm{Con}\,(\mathbf A,\tau)$.
\end{proof}

We remind the following Mal'cev Theorem, \cite[Lem 3.1]{BuSa}.

\begin{theorem}
Let $\mathbf A$ be an algebra and $\phi\subseteq A^2$. Then
$(a,b)\in\Theta(\phi)$ if and only if there exist two finite
sequences of terms $t_1(\overline x_1,x), \ldots, t_n(\overline
x_n,x)$ and pairs $(a_1,b_1),\ldots,(a_n,b_n)\in\phi$ with
$$
a=t_1(\overline x_1,a_1),\, t_i(\overline x_i,b_i)=t_{i+1}(\overline
x_{i+1},a_{i+1}) \mbox{ and } t_n(\overline x_n,b_n)=b
$$
for some $\overline x_1,\ldots,\overline x_n\in A$.
\end{theorem}

We say that an algebra $\mathbf B$ has the  Congruence Extension
Property (CEP for short) if, for any algebra $\mathbf A$ such that
$\mathbf B$ is a subalgebra of $\mathbf A$ and for any congruence
$\theta \in \mathrm{Con}\, \mathbf B,$ there is a congruence $\phi
\in \mathrm{Con}\, \mathbf A$ such that $\theta = (B \times B)\cap
\phi.$ A variety $\mathcal K$ has the CEP if every algebra in
$\mathcal K$ has the CEP.  For example, the variety of MV-algebra,
or the variety of BL-algebras or the variety of state-morphism
MV-algebras (see \cite[Lem 6.1]{DKM}) satisfies the CEP.

\begin{theorem}\label{cep}
A variety $\mathcal V_\tau$ satisfy the CEP if and only if $\mathcal
V$ satisfies the CEP.
\end{theorem}

\begin{proof}
Let us have a variety $\mathcal V$ with the CEP. If $\mathbf
A\in\mathcal V$ is such that $(\mathbf A,\tau)$ is an algebra with
state-morphism,  for any subalgebra $(\mathbf B,\tau)\subseteq
(\mathbf A,\tau)$ and any $\phi\in\mathrm{Con}\,(\mathbf B,\tau),$
the condition $\phi=B^2\cap \Theta(\phi)$ holds.

Now we prove $\Theta(\phi)=\Theta_\tau(\phi)$. To show that, assume
$(a,b)\in \Theta(\phi).$   Mal'cev's Theorem shows the existence of
finite sequences of terms $t_1(\overline x_1,x), \ldots,
t_n(\overline x_n,x)$ and pairs $(a_1,b_1),\ldots,(a_n,b_n)\in\phi$
with
$$
a=t_1(\overline x_1,a_1),\, t_i(\overline x_i,b_i)=t_{i+1}(\overline
x_{i+1},a_{i+1}) \mbox{ and } t_n(\overline x_n,b_n)=b
$$
for some $\overline x_1,\ldots,\overline x_n\in A$. Because $\tau$
is an endomorphism, we obtain also equalities
$$
\tau (a)=t_1(\tau (\overline x_1),\tau (a_1)),\, t_i(\tau (\overline
x_i),\tau(b_i))=t_{i+1}(\tau(\overline x_{i+1}),\tau(a_{i+1}))$$ and
$$
t_n(\tau (\overline x_n),\tau (b_n))=\tau (b).
$$
We have assumed that $\phi\in\mathrm{Con}\,(\mathbf B,\tau)$, thus
$(a_i,b_i)\in\phi$ yields $(\tau (a_i),\tau (b_i))\in\phi$ for any
$i=1,\ldots,n.$ Now, we have obtained $(\tau(a),\tau(b))\in\Theta
(\phi)$. In other words, $\Theta (\phi)\in\mathrm{Con}\,(\mathbf
A,\tau)$ and thus $\Theta (\phi)=\Theta_\tau (\phi)$.

If $\mathcal V_\tau$ has the CEP, then for any $\mathbf A\in\mathcal
V,$ we have $\mathrm{Con}\,\mathbf A = \mathrm{Con}\,(\mathbf
A,\Id_A)$. Clearly, the CEP on $(\mathbf A,\Id_A)$  yields the CEP
on $\mathbf A$.
\end{proof}

\section{Applications to Special Types of Algebras}

In this section, we apply  a general result concerning generators of
some varieties of state-morphism algebras, Theorem \ref{th:3.3.1},
to the variety of state-morphism BL-algebras, state-morphism
MTL-algebras, state-morphism non-associative BL-algebras, and
state-morphism pseudo MV-algebras, when we use different systems of
t-norms on the real interval $[0,1]$ and a special type of pseudo
MV-algebras, respectively.

Algebras for which the logic MTL is sound are called MTL-algebras.
They can be characterized as prelinear commutative bounded integral
residuated lattices. In more detail, according to \cite{EsGo}, an
algebraic structure $\mathbf A=(A;\wedge,\vee,\ast,\rightarrow,0,1)$
of type $\langle 2,2,2,2,0,0\rangle$ is an {\it MTL-algebra} if

\begin{itemize}
\item[(M1)] $(A;\wedge,\vee,0,1)$ is a bounded lattice with the top element 0
and bottom element 1,
\item[(M2)] $(A;\ast,1)$ is a commutative monoid,
\item[(M3)] $\ast$ and $\rightarrow$ form an adjoint pair, that is, $z*x\le y$
if and only if $z\le x\rightarrow y,$ where $\le$ is the lattice
order of $(A;\wedge,\vee)$ for all $x, y, z \in A,$ (the residuation
condition),
\item[(M4)] $(x\rightarrow y)\vee(y\rightarrow x)=1$ holds for all $x,y \in A$
(the prelinearity condition).
\end{itemize}

If $t$ is any left-continuous t-norm on $[0,1],$ we define two
binary operations  $\ast_t$  $\to_t$ on $[0,1]$ via   $x\ast_t
y=t(x,y)$ and $x\to_t y = \sup\{z\in [0,1]: t(z,x)\le y\}$ for $x,y
\in [0,1],$ then $\mathbb I_t =([0,1]; \min, \max,
\ast_t,\to_t,0,1)$ is an example of an MTL-algebra. An MTL-algebra
$\mathbb I_t$ is a BL-algebra iff $t$ is continuous.

Due to \cite{EsGo},  the class $\mathcal T_{lc},$ which denotes the
system of all BL-algebras $\mathbb I_t,$ where $t$ is a
left-continuous t-norm on the interval $[0,1],$ generates the
variety of MTL-algebras. This result was strengthened in \cite{Vet}
who introduced the class of regular left-continuous t-norms which is
strictly smaller than the class of left-continuous t-norms, but they
generate  the variety of MTL-algebras.

According to \cite{Bot}, we say that an algebra $\mathbf
A=(A;\vee,\wedge,\cdot,\rightarrow,0,1)$ of type $\langle
2,2,2,2,0,0\rangle$ is a {\it non-associative BL-algebra}
(naBL-algebra in short) if
\begin{itemize}
\item[(A1)] $(A;\vee,\wedge,0,1)$ is a bounded lattice,
\item[(A2)] $(A;\cdot,1)$ is a commutative groupoid with the neutral element
$1$,
\item[(A3)] any $x,y,z\in A$ satisfy $x\cdot y\leq z$ if and only if $x\leq
y\rightarrow z,$
\item[(A4)] algebra satisfy the divisibility axiom ($x\cdot (x\rightarrow
y)=x\wedge y$),
\item[(A5)] algebra satisfy the $\alpha$-prelinearity and $\beta$-prelinearity
($x\rightarrow y\vee \alpha^a_b(y\rightarrow x)=x\rightarrow y\vee
\beta^a_b(y\rightarrow x)=1$), where $\alpha^a_b(x)=(a\cdot
b)\rightarrow(a\cdot (b\cdot x))$ and $\beta^a_b(x)=b\rightarrow
(a\rightarrow ((a\cdot b)\cdot x))$.
\end{itemize}

A function $t: [0,1]\times [0,1] \to [0,1] $ on the  interval
$[0,1]$ of reals is said to be a {\it non-associative} t-norm
(nat-norm briefly) if
\begin{itemize}
\item[(nat1)] $([0,1];t,1)$ is a commutative groupoid with the neutral element
$1$,
\item[(nat2)]  $t$ is continuous in the usual sense,
\item[(nat3)] if $x,y,z\in [0,1]$ are such that $x\leq y,$ then $t(x,z)\leq
t(y,z)$.
\end{itemize}

According to \cite[Thm 5]{Bot}, for any nat-norm there is a unique
binary operation $\rightarrow_t$ satisfying the adjointness
condition, i.e. $t(x,y)\leq z$ if and only if $x\leq y\rightarrow_t
z$. Moreover, an algebra $\mathbb
I_t^{na}:=([0,1];\min,\max,t,\rightarrow_t,0,1)$ is an naBL-algebra.

The class of all naBL-algebras is denoted by $na\mathcal{BL}$ and
$na\mathcal T$ denotes the class of all naBL-algebras $\mathbb
I^{na}_t$ for any non-associative t-norm. The main result on
non-associative BL-algebras says that $na\mathcal T$ is the
generating class for the variety $na\mathcal{BL},$ \cite[Thm
8]{Bot}:

\begin{theorem}\label{m:5.1} There hods
$$
na\mathcal{BL} = \mathsf{IP_SSP_U}(na\mathcal T).
$$
\end{theorem}

Finally, we recall that a noncommutative generalization of
MV-algebras was introduced in \cite{GeIo} as {\it pseudo
MV-algebras} or in \cite{Rac} as {\it generalized MV-algebras}.
According to \cite{Dvu}, every pseudo MV-algebra $(M;\oplus,
^-,^\sim,0,1)$ of type $\langle 2,1,1,0,0\rangle$ is an interval in
a unital $\ell$-group $(G,u)$ with strong unit $u,$ i.e. $M\cong
\Gamma(G,u):=[0,u],$ where $x\oplus y =(x+y)\wedge,$ $x^-=u-x,$
$x^\sim =-x+u,$ $0=0,$ and $1=u.$ If $(G,u)$ is double transitive
(for definitions and details see \cite{DvHo}), then $\Gamma(G,u)$
generates the variety of pseudo MV-algebras, \cite[Thm 4.8]{DvHo}.
For example, if $\mathrm{Aut}(\mathbb R)$ is the set of all
automorphisms of the real line $\mathbb R$ preserving the natural
order in $\mathbb R$ and $u(t):=t+1,$ $t \in \mathbb R$,  let
$\mathrm{Aut}_u(\mathbb R)=\{ g \in \mathrm{Aut}(\mathbb R): g \le
nu$ for some integer $n\ge 1 \}.$ Then
$\Gamma(\mathrm{Aut}_u(\mathbb R),u)$ is double transitive and it
generates the variety of pseudo MV-algebras, see \cite[Ex
5.3]{DvHo}.

Now we apply the general statement, Theorem \ref{th:3.4.1}, on
generators to different types of state-morphism algebras. We recall
that $\mathcal T$ was defined as the class of all BL-algebras
$\mathbb I_t,$ where $t$ is a continuous t-norm on $[0,1].$

\begin{theorem}\label{ad:5.2} {\rm (1)}  The variety of all
state-morphism MV-algebras is generated by the diagonal
state-morphism MV-algebra $D([0,1]_{MV}).$

{\rm (2)} The variety of all state-morphism BL-algebras is generated
by the class $\{D(\mathbb I_t): \mathbb I_t \in\mathcal T\}.$

{\rm (3)}  The variety of all state-morphism MTL-algebras is
generated by the class $\{D(\mathbb I_t): \mathbb I_t \in\mathcal
{T}_{lc}\}.$

{\rm (4)} The variety of all state-morphism naBL-algebras is
generated by the class $\{D(\mathbb I^{na}_t): \mathbb I_t\in
na\mathcal T\}.$

{\rm (5)}  If a unital $\ell$-group $(G,u)$ is double transitive,
then $D(\Gamma(G,u))$ generates the variety of  state-morphism
pseudo MV-algebras.
\end{theorem}

\begin{proof}
(1) It follows from the fact that the MV-algebra of the real
interval $[0,1]$ generates the variety of MV-algebras, see e.g.
\cite[Prop 8.1.1]{CDM}, and then apply Theorem \ref{th:3.4.1}.

(2) The statement follows from the fact that $\mathsf{V}(\mathcal
T)$ is by \cite[Thm 5.2]{CEGT} the variety $\mathcal {BL}$ of all
BL-algebras. Now it suffices to apply Theorem \ref{th:3.4.1}.

(3) By \cite{EsGo}, the class $\mathcal T_{lc}$ of all $\mathbb
I_t,$ where $t$ is any left-continuous t-norms on the interval
$[0,1],$ generates the variety of MTL-algebras; then apply Theorem
\ref{th:3.4.1}.

(4) By \cite[Thm 8]{Bot} or Theorem \ref{m:5.1}, the class
$na\mathcal{T}$ of all $\mathbb I_t,$ where $t$ is any
non-associative t-norms on the interval $[0,1],$ generates the
variety of non-associative BL-algebras; then apply again Theorem
\ref{th:3.4.1}.

(5) By the above, $\Gamma(G,u)$ generates the variety of pseudo
MV-algebras, see also \cite[Thm 4.8]{DvHo}; then apply Theorem
\ref{th:3.4.1}.
\end{proof}

We note that the case (1) in Theorem \ref{th:3.4.1} was an open
problem posed in \cite{DiDv} and was positively solved in \cite[Thm
5.4(3)]{DKM}.

\section{Conclusion}

In the paper, we have presented a general approach to theory of
state-morphism algebras which generalizes state-morphism MV-algebras
and state-morphism BL-algebras as pairs $(\mathbf A,\tau),$ where
$\mathbf A$ is an algebra of type $F$ and $\tau$ is an endomorphism
of $\mathbf A$ such that $\tau \circ \tau = \tau.$

This enables us to present  complete characterizations of
subdirectly irreducible state BL-algebras and subdirectly
irreducible state-morphism BL-algebras, Theorem \ref{th:5.1}, which
generalizes the results from \cite{DiDv, DDL2, Dvu1, DKM}.

A general approach is studied in the third section where the main
result, Theorem \ref{diag}, says that every subdirectly irreducible
state-morphism algebra can be embedded into a diagonal one.

The fourth section describes some generators of the varieties of
state-morphism algebras, and Theorem \ref{th:3.4.1} shows that if a
class $\mathcal K$ generates a variety $\mathcal V$ of algebras of
the same type $F$, then the variety of state-morphism algebras whose
$F$-reduct belongs to the class $\mathcal K$ is generated by the
class of diagonal state-morphism algebras $D(\mathbf A),$ where
$\mathbf A \in \mathcal K.$  In addition, Theorem \ref{cep} deals
with the CEP for the variety of state-morphism algebras.

In Theorem \ref{ad:5.2}, Theorem \ref{th:3.4.1} was applied to the
special class of algebras: MV-algebras, BL-algebras, MTL-algebras,
non-associative BL-algebras, and pseudo MV-algebras to obtain the
generators of the corresponding varieties of state-morphism
algebras.

During the study on this paper, we found some interesting open
problems like: (1) find a characterization of an analogue of a
state-operator that is not necessarily a state-morphism operator,
(2) if the lattice of varieties of some variety is countable, how
big is the lattice of corresponding state-morphism algebras, e.g. in
the case of MV-algebras, the lattice under question is uncountable
\cite{DKM}, (3) decidability of the variety of state-morphism
algebras.

\end{document}